\documentclass[11pt]{amsart}
\usepackage[obeyclassoptions=true,subpreambles=false]{standalone}
\usepackage{import}

\usepackage[utf8]{inputenc}
\usepackage[T1]{fontenc}
\usepackage{lmodern}
\usepackage{mathtools,amsfonts,amsthm,mathrsfs,amssymb,bm}
\usepackage[hmargin=1in,letterpaper]{geometry}

\usepackage{subcaption}
\usepackage{float}

\usepackage{textcomp}
\mathtoolsset{showonlyrefs}

\usepackage[shortlabels]{enumitem}

\usepackage[color]{showkeys}

\PassOptionsToPackage{pdfusetitle}{hyperref}
\usepackage{bookmark}
\usepackage{hypcap}
\usepackage{microtype}

\newtheorem{theorem}{Theorem}

\newtheorem{conjecture}[theorem]{Conjecture}
\newtheorem{question}[theorem]{Question}

\theoremstyle{definition}
\newtheorem{definition}{Definition}

\newcommand*{\lam}{\lambda}

\renewcommand{\epsilon}{\varepsilon}
\newcommand{\eps}{\varepsilon}

\renewcommand{\subset}{\subseteq}

\newlist{enumerata}{enumerate}{1}
\setlist[enumerata]{label=\upshape{(\alph*)}}
\setlist[enumerate]{label=\upshape{(\roman*)}}

\title{On the occupancy fraction of the antiferromagnetic Ising model}
\date{\today}

\author[E.\ Davies]{Ewan Davies}
\thanks{Supported in part by NSF grant CCF-2309707}
\address{Department of Computer Science, Colorado State University, Fort Collins, USA}
\email{research@ewandavies.org}

\author[O.\ LeBlanc]{Olivia LeBlanc}
\address{Department of Computer Science, Colorado State University, Fort Collins, USA}
\email{Olivia.Leblanc@colostate.edu}

\begin{document}
\begin{abstract}
We study the maximum and minimum occupancy fraction of the antiferromagnetic Ising model in regular graphs. The minimizing problem is known to determine a computational threshold in the complexity of approximately sampling from the Ising model at a given magnetization, and our results determine this threshold for nearly the entire relevant parameter range in the case $\Delta=3$. A small part of the parameter range lies outside the reach of our methods, and it seems challenging to extend our techniques to larger $\Delta$. 
\end{abstract}

\maketitle

\section{Introduction}

In computer science, the Ising model is an important probability distribution over cuts in a graph which models magnetic material.
One of the canonical computational problems associated with the Ising model is to approximate the partition function, which is a weighted sum over cuts in a graph. 
Naive computation of this quantity takes exponential time, and we are interested in obtaining an approximate answer in polynomial time. 
Exact computation is also known to be \#P-hard.
It is natural to impose a global constraint on the model which fixes the magnetization, or equivalently the sizes of each side of the cut. 
Depending on the parameters of the model, it has been established~\cite{DP23} that there is a \emph{computational threshold} in the complexity of the approximate counting problem at given magnetization. 
That is, there is a boundary in the parameter space of the problem such that on one side of the boundary there is an efficient algorithm, while unless P=NP there cannot be such an algorithm on the other side of the boundary. 

The study of computational thresholds in globally-constrained approximate counting problems was initiated in~\cite{DP23}, where a threshold was established for the problem of approximately counting independent sets of a given density in bounded-degree graphs. 
The boundary in the parameter space is a simple one: there is a single density parameter $\alpha$ and for each maximum degree $\Delta\ge 3$ there is an explicit critical density $\alpha_c(\Delta)$ such that various efficient algorithmic approaches prevail for $\alpha<\alpha_c$ (see~\cite{DP23,JMPV23,JPSS22}) and unless some complexity-theoretic collapse occurs no efficient algorithms prevail when $\alpha>\alpha_c$.

The methods of~\cite{DP23} extend to the antiferromagnetic Ising model and show for the class of $\Delta$-regular graphs that there is a computational threshold of the type outlined above, though rather curiously the method does not determine the location of the threshold explicitly. 
Our main contribution is to determine the location of this threshold in the case of $3$-regular graphs when the edge activity of the antiferromagnetic Ising model is strong enough. 
In principle our techniques could extend to larger degrees, though the method seems to require prohibitively expensive computation in these cases. Our method seems to break down for more fundamental reasons as the edge activity weakens.

We define the Ising model precisely before continuing the discussion.
For an edge activity parameter $B$ and external field $\lambda$, the partition function $Z_G(B,\lambda)$ of the Ising model on a graph $G$ is given by
\[ Z_G(B,\lambda) := \sum_{\sigma : V(G)\to\{+,-\}}B^{m_G(\sigma)}\lambda^{n_G(+, \sigma)}, \]
where $\sigma$ is an assignment of a spin in $\{+,-\}$ to each vertex of $G$, $m_G(\sigma)$ is the number of edges of $G$ whose endpoints get the same spin under $\sigma$ (we call such edges \emph{monochromatic}) and $n_G(+,\sigma)$ is the number of vertices which receive the spin $+$ under $\sigma$.
The model is antiferromagnetic if $B\in(0,1)$, and without loss of generality we consider $0<\lambda\le 1$ because the model is symmetric under swapping the spins and taking $\lambda\mapsto 1/\lambda$.
We define the Ising measure $\mu_{G,B,\lambda}$ on spin assignments $\sigma : V(G)\to \{+,-\}$ via 
\[ \mu_{G,B,\lambda}(\sigma) = \frac{B^{m_G(\sigma)}\lambda^{n_G(+,\sigma)}}{Z_G(B,\lambda)}. \]
The Ising measure can be defined on the infinite $\Delta$-regular tree by the DLR equations, and in general for fixed $B$ and $\lambda$ one can obtain multiple measures this way (see e.g.~\cite{FV17} for an introduction to these ideas in the case of the ferromagnetic Ising model and other spin systems). 
There is a region, known as the \emph{Gibbs uniqueness} region in the $(B,\lambda)$ parameter space such that the resulting measure is unique. 
For $\Delta\le 2$ the measure is always unique, but for $\Delta\ge 3$ some interesting behavior arises. 
Let $B_c(\Delta)=(\Delta-2)/\Delta$. Then for $B\in (B_c, 1)$ the measure is unique for all values of $\lambda$, but for $B\in(0,B_c)$ there is a critical external field $\lambda_c(\Delta,B)$ such that for $\lambda\in(0,\lambda_c]$ the measure is unique and for $\lambda \in (\lambda_c,1]$ multiple measures are possible.
We record the explicit function giving $\lambda_c$ here:
\[ \lambda_c(\Delta,B) = \frac{1-\sqrt{r/s}}{1+\sqrt{r/s}}\left(\frac{1+\sqrt{rs}}{1-\sqrt{rs}}\right)^{\frac{B_c-1}{B_c+1}},\]
where $B_c=B_C(\Delta)$ and 
\begin{align*}
    r &= \frac{B_c-B}{B_c+B},&
    s &= \frac{1-B}{1+B}.
\end{align*}
The fact that this is not an especially simple function of $B$ adds to the difficulty of the problems we study.

The magnetization of a spin assignment $\sigma$ is $M(\sigma)=\sum_{v\in V}\sigma_v$, and given a fixed graph, fixing the magnetization is equivalent to fixing $n_G(+,\sigma)$. 
If we view the partition function as a polynomial in $\lambda$ we have 
\[ Z_G(B,\lambda) = \sum_{k=0}^{|V|} c_k \lambda^k \]
for some coefficients $c_k = c_k(G,B)$ which depend on the graph and the edge activity parameter $B$. 
We write $\mathcal{B}_k(G)$ for the subset of spin assignments $\sigma$ such that $n_G(+,\sigma)$ is exactly $k$, giving 
\[ c_k = \sum_{\sigma\in\mathcal{B}_k(G)}B^{m_G(\sigma)}. \]

For each given $B$, we are interested in establishing a threshold in $k$ for the computational problem of approximating $c_k(G,B)$ in the class of $\Delta$-regular graphs. 
To discuss the threshold it is convenient to normalize by the number of vertices in the graph and define the \emph{occupancy fraction}
\[ \alpha_G(B,\lambda) = \frac{\lambda}{|V(G)|} \frac{\partial}{\partial \lambda} \log Z_G(B,\lambda). \]
It is a straightforward calculation to show that this is the expected fraction of the vertices which receive spin $+$ in a sample from $\mu_{G,B,\lambda}$.
The methods of~\cite{DP23} establish in some generality that for two-spin antiferromagnetic systems, minimizing $\alpha$ over a class of graphs of interest determines the computational threshold. 
Let $\mathcal{G}_\Delta$ be the class of $\Delta$-regular graphs and let
\[ \alpha_{\inf}(\Delta, B,\lambda) = \inf_{G\in \mathcal{G}_\Delta} \alpha_G(B,\lambda). \]
The result of~\cite{DP23} relevant to the antiferromagnetic Ising model in regular graphs establishes a computational threshold $\alpha_c$ as follows. 
Note that a \emph{fully polynomial-time ramdomized approximation scheme} or FPRAS for a quantity $Q(x)$ is a family of randomized algorithms for each $\eps>0$ that compute with probability at least 3/4 a value $\hat Q(x)$ such that $e^{-\eps} \le \hat Q(x) /Q(x)\le e^{\eps}$, in time polynomial in the input size (e.g. the number of vertices of the graph represented by $x$) and $1/\eps$.

\begin{theorem}[{Davies and Perkins~\cite[Thm.\ 3]{DP23}}]
For $\Delta\ge 3$ and $0<B<B_c(\Delta)$, we write $\alpha_c = \alpha_{\inf}(\Delta, B,\lambda_c(\Delta,B))$. Then $\alpha_c$ is a computational threshold in the following sense.
\begin{enumerate}[{label=\textup{(\alph*)}}]
    \item For every $\alpha<\alpha_c$ there is an FPRAS for $c_{\lfloor \alpha n\rfloor}(G,B)$ for all $n$-vertex $\Delta$-regular graphs $G$. 
    \item Unless NP=RP, for every $\alpha\in(\alpha_c,1/2]$ there is no FPRAS for $c_{\lfloor \alpha n\rfloor}(G,B)$ for $n$-vertex $\Delta$-regular graphs $G$.
\end{enumerate}
\end{theorem}

Given this theorem, to determine the computational threshold requires finding an explicit formula for $\alpha_c$. 
In~\cite{DP23} it was conjectured that for $\Delta\ge 3$ and $B\in (0,B_c(\Delta))$, $\alpha_c$ is the occupancy fraction of the complete graph $K_{\Delta+1}$ evaluated at edge activity $B$ and external field $\lambda_c(\Delta,B)$. We confirm this conjecture in the case $\Delta=3$ and $B\in(0,0.3128)$, noting that in the case of $\Delta=3$ the conjecture covers values of $B$ up to $B_c(3)=1/3$.

\begin{theorem}\label{thm:d3lc}
    The complete graph $K_{4}$ minimizes the occupancy fraction $\alpha_G(B,\lambda_c(3,B))$ over $3$-regular graphs when $0\le B\le 0.3128$.
\end{theorem}

This is an example of an extremal problem in combinatorics: we seek the minimum of some graph-theoretic quantity (the occupancy fraction) over a class of graphs. 
Such problems are well-studied in general, though to our knowledge the case of the antiferromagnetic Ising model has escaped significant attention. 
Inspired by this wider interest in graph theory, we also tackle the minimization question in $3$-regular graphs over the entire parameter range $(B,\lam)\in[0,1]^2$. 
While we obtain partial results, we gather evidence that the problem is significantly more challenging than analogous problems involving other well-known spin systems.
Earlier results on extremal problems involving the antiferromagnetic Ising model on regular graphs typically focus on a different parameter, either the partition function itself or equivalently the \emph{free energy density} given by
\[ F_G(B,\lambda) = \frac{1}{|V(G)|}\log Z_G(\lambda). \]
The free energy density includes a convenient normalization factor that facilitates the comparison of graphs on different numbers of vertices (cf.\ the same normalization in the definition of the occupancy fraction).
It is easy to see that 
\[ F_G(B,\lambda) = \int_{0}^{\lambda}\frac{1}{\ell}\alpha_G(B,\ell) \,\mathrm d\ell, \]
and hence an extremal result for the occupancy fraction in a suitably `downward-closed' region $\{(B,\ell) : 0\le \lam\le \ell\}$ immediately implies an extremal result for the free energy density.

\begin{theorem}\label{thm:d3min}
    The complete graph $K_4$ minimizes the occupancy fraction $\alpha_G(B,\lambda)$ over $3$-regular graphs in the union of the regions below:
    \begin{enumerate}
        \item $\mathcal{R}^{\min}_1 := \{ (B,\lambda) : \frac{59}{100} \le B \le 1 \land0 \le \lambda \le  \frac{177}{200}B-\frac{267}{1000}B^2\}$
        \item $\mathcal{R}^{\min}_2 := \{ (B,\lambda) : 0 \le B \le \frac{59}{100} \land 0 \le \lambda \le \frac{9}{10}B\}$
        \item $\mathcal{R}^{\min}_3 := \{ (B,\lambda) : \frac{1}{10} \le B \le \frac{1}{2} \land \frac{1}{2}B \le \lambda \le \frac{99}{100}B\}$
        \item $\mathcal{R}^{\min}_4 := \{ (B,\lambda) : \frac{1}{2} \le B \le \frac{9}{10} \land \frac{3}{20} \le \lambda \le \min\big\{B\frac{53}{100}\big\}$
        \item $\mathcal{R}^{\min}_5 := \{ (B,\lambda) : \frac{3}{5} \le B \le \frac{49}{50} \land \frac{7}{20} + \frac{11}{50}B \le \lambda \le \frac{3}{8}+\frac{6}{25}B\}$
        \item $\mathcal{R}^{\min}_6 := \{ (B,\lambda) : \frac{1}{4} \le B \le \frac{9}{20} \land \frac{49}{50}B \le \lambda \le -\frac{17}{200}+\frac{7}{5}B\}$
        \item $\mathcal{R}^{\min}_7 := \{ (B,\lambda) : \frac{56}{125} \le B \le \frac{3}{5} \land \frac{99}{100}B \le \lambda \le \frac{133}{200} - \frac{13}{50}B\}$
        \item $\mathcal{R}^{\min}_8 := \{ (B,\lambda) : \frac{13}{50} \le B \le \frac{9}{25} \land -\frac{3}{40}+\frac{34}{25}B \le \lambda \le -\frac{9}{100}+\frac{36}{25}B\}$
    \end{enumerate}

    The complete graph $K_4$ also minimizes the free energy density $F_G(B,\lambda)$ over $3$-regular graphs in the union of these regions.

    The complete graph $K_4$ does not minimize the occupancy fraction $\alpha_G(B,\lambda)$ in the entire region $[0,1]^2$. There is a region of the parameter space $\mathcal{R}^{\mathrm{Pet}}$ where the Petersen graph has smaller occupancy fraction. 
    It is also false that the minimizer is always either $K_4$ or the Petersen graph, as a third graph which we call the Goose graph\footnote{The name arises from a graph6 representation of the graph: \texttt{I\}GOOSE@W}.} has smaller occupancy fraction than either of these two in a region $\mathcal{R}^{\mathrm{Goose}}$ of $[0,1]^2$.
\end{theorem}

\begin{figure}[htp]
    \centering
    \vspace*{2ex}
    \includegraphics[width=0.6\textwidth]{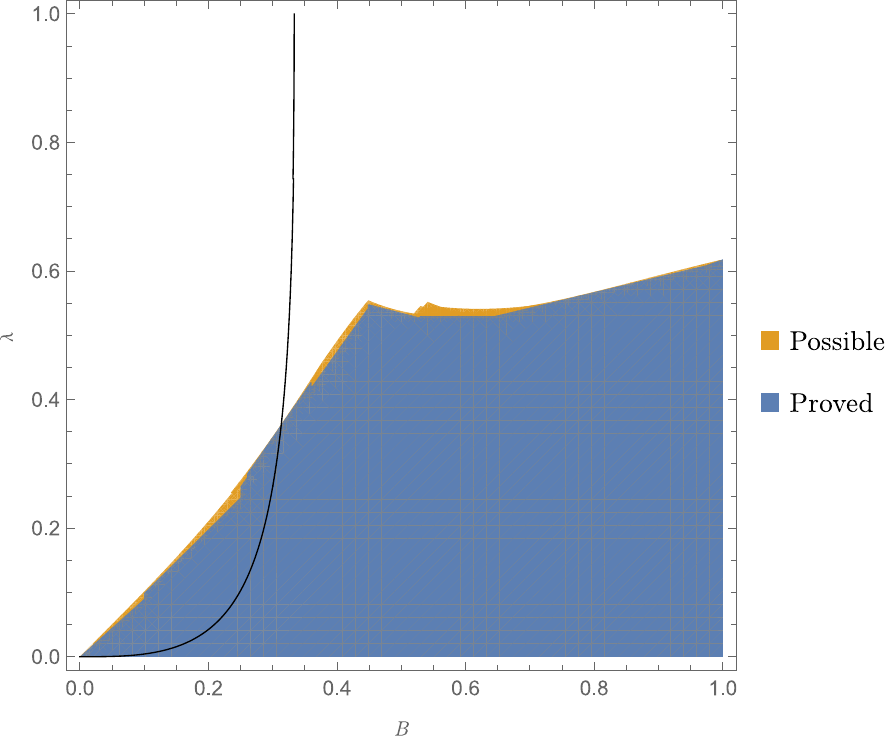}
    \vspace*{2ex}
    \caption{Minimizing the occupancy fraction over 3-regular graphs.
    The blue region is union of the regions in Theorem~\ref{thm:d3min} where we proved that $K_4$ minimizes the occupancy fraction. The yellow region is where we numerically evaluated that it is theoretically possible for our method to show that $K_4$ minimizes the occupancy fraction, but we did not push the symbolic proofs of dual feasibility to the limit.
    The black line is $\lam_c(3,B)$.}
    \label{fig:d3min}
\end{figure}

\begin{figure}[tp]
    \centering
    \vspace*{8ex}
    \includegraphics[width=0.6\linewidth]{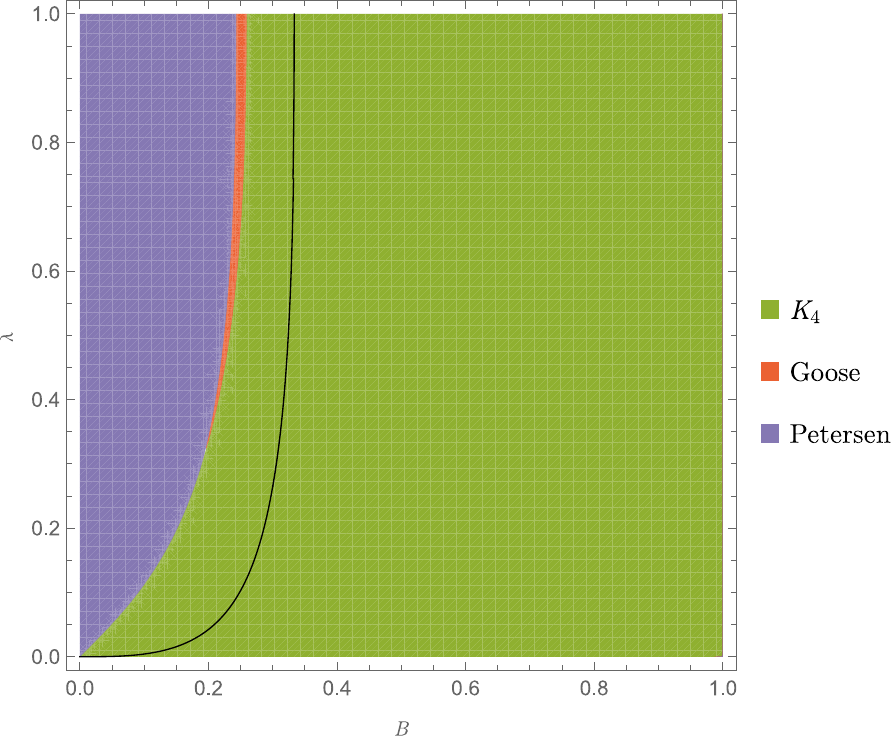}
    \\\vspace*{8ex}
    \includegraphics[width=0.2\textwidth]{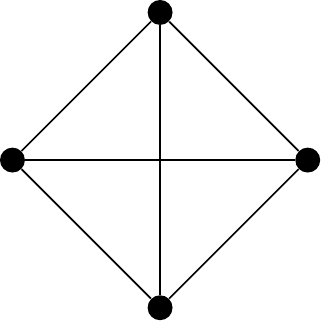}
    \hspace*{0.1\textwidth}
    \includegraphics[width=0.2\textwidth]{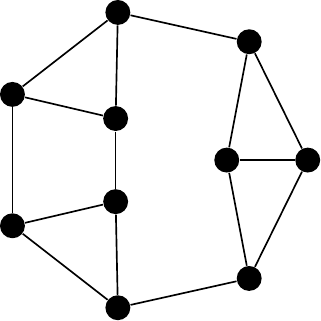}
    \hspace*{0.1\textwidth}
    \includegraphics[width=0.2\textwidth]{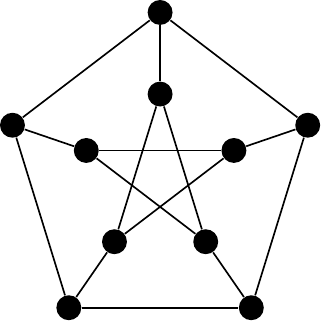}
    \\\vspace*{8ex}
    \caption{Considering only the pictured graphs, $K_4$, the Goose graph, and the Petersen graph respectively, the plot shows the region where each of them has the smallest occupancy fraction out of the three. 
    The black line is $\lam_c(3,B)$. 
    This shows that it is not possible to strengthen the relevant part of Theorem~\ref{thm:d3min} by enlarging the regions to cover the entire parameter range $[0,1]^2$, but it does not disprove the original conjecture in~\cite{DP23} about $K_4$ minimizing when $\lambda=\lambda_c(3,B)$.}
    \label{fig:d3posmin}
\end{figure}

We show the union of the regions in Theorem~\ref{thm:d3min} in Figure~\ref{fig:d3min}, and a plot of where the Petersen and Goose graphs have smaller occupancy fraction in Figure~\ref{fig:d3posmin}.

Our methods yield insights into the corresponding problem of maximizing the occupancy fraction of the antiferromagnetic Ising model too.

\begin{theorem}\label{thm:d3max}
    The complete bipartite graph $K_{3,3}$ maximizes the occupancy fraction $\alpha_G(B,\lambda)$ over $3$-regular graphs in the union of the regions below:
    \begin{enumerate}
        \item $\mathcal{R}^{\max}_1 := \{ (B,\lambda) : 0 \le B \le 1/5 \land 0 \le \lambda \le 3B/10\}$,
        \item $\mathcal{R}^{\max}_2 := \{ (B,\lambda) : 1/5 \le B \le 2/5 \land 0 \le \lambda \le 4B/10\}$,
        \item $\mathcal{R}^{\max}_3 := \{ (B,\lambda) : 2/5 \le B \le 1 \land 0 \le \lambda \le 5B/10\}$,
    \end{enumerate}
\end{theorem}

\begin{figure}[tp]
    \centering
    \vspace*{2ex}
    \includegraphics[width=0.6\textwidth]{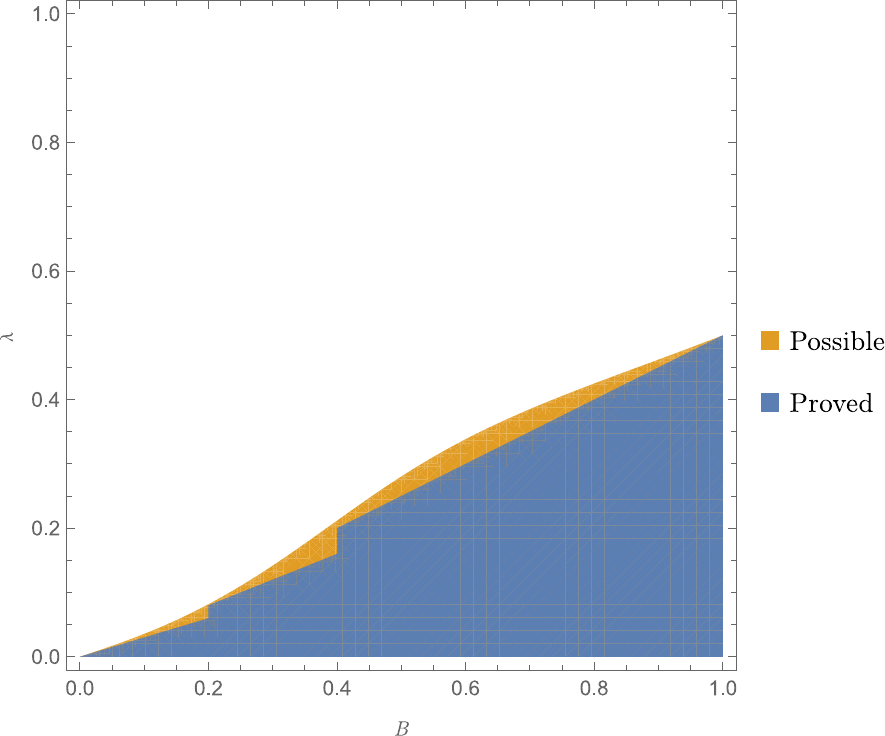}
    \vspace*{2ex}
    \caption{Maximizing the occupancy fraction over 3-regular graphs.
    The blue region is union of the regions in Theorem~\ref{thm:d3max} where we proved that $K_{3,3}$ maximizes the occupancy fraction. The yellow region is where we numerically evaluated that it is theoretically possible for our method to show that $K_{3,3}$ maximizes the occupancy fraction, but we did not push the symbolic proofs of dual feasibility to the limit.}
    \label{fig:d3max}
\end{figure}

We show the union of the regions in Theorem~\ref{thm:d3max} in Figure~\ref{fig:d3max}.
We conjecture that the above result extends to all $\Delta\ge 3$ and the entire parameter range of the antiferromagnetic Ising model.
Despite proving that $K_4$ cannot minimize the occupancy fraction on the entire parameter space, we also conjecture that the complete graph minimizes the free energy density.
Given the complexities apparent in the 3-regular case, these may be considered somewhat bold conjectures.

\begin{conjecture}\label{conj:d3}
    For all $\Delta\ge 3$ and all $(B,\lambda)\in[0,1]^2$, the following hold.
    \begin{enumerate}
        \item\label{itm:conjd3max} The complete bipartite graph $K_{\Delta,\Delta}$ maximizes the occupancy fraction $\alpha_G(B,\lambda)$ over $\Delta$-regular graphs.
        \item\label{itm:conjd3min} The complete graph $K_{\Delta+1}$ minimizes the free energy density $F_G(B,\lambda)$ over $\Delta$-regular graphs.
    \end{enumerate}
\end{conjecture}

\subsection{Related work}

An analogous algorithmic problem for the ferromagnetic Ising model was solved in~\cite{CDKP22}. 
Therein, the authors showed that there is a computational threshold denoted $\eta_c$ in the mean magnetization for the problem of approximately counting and sampling from the model at fixed magnetization. 
There is an analogous extremal graph theory problem to determine $\eta_c$ which corresponds to maximizing the mean magnetization over bounded-degree graphs, and in~\cite{CDKP22} the solution was shown to be the `zero-field + measure' on the infinite regular tree. 
The ferromagnetic Ising model has rather different behavior as there is an FPRAS for the partition function for all values of the parameters and so one cannot use a reduction from the unconstrained model to the model at given magnetization to show hardness. 
The main result of~\cite{CDKP22} was to demonstrate that nonetheless, constraining the magnetization results in a hard computational problem in a natural parameter range. 
The relevant extremal problem is related to the one here of minimizing the occupancy fraction, but for the ferromagnetic Ising model and over the class of bounded-degree graphs. 
For the ferromagnetic model there are correlation inequalities known as the GKS inequalities that can be harnessed in the proof, and we are not aware of analogous techniques that apply to the antiferromagnetic Ising model. 

The algorithmic problem of approximating individual coefficients of a (grand canonical) partition function an be approached in various ways. 
Markov chain methods for approximate sampling are the oldest approach~\cite{BD97}, though the traditional path coupling technique does not always give optimal bounds. 
The main algorithmic idea in~\cite{DP23} is to perform rejection sampling, using the fact that the constrained model is simply the unconstrained model conditioned on the constraint. 
The key to the analysis of such an algorithm is the central limit theorem provided by a zero-free region of the partition function.
Later works demonstrated the viability of deterministic and randomized algorithms inspired by \emph{local} central limit theorems~\cite{JPSS22}, and subsequent works showed that more advanced techniques can improve the analysis of local Markov chains for such problems~\cite{AL20a,JMPV23,KPPY24}.
Rather surprisingly, there are advanced, general techniques for relating algorithms for approximate counting in constrained and unconstrained spin systems~\cite{HK24}. 
It is not clear that such methods are sharp enough to witness computational thresholds of the type we study here, but given knowledge of relevant computational thresholds (and a few additional related properties), these general approaches often yield black-box algorithms with good performance despite having relatively little problem-specific structure.

Extremal problems of the type we address here are also well-studied. 
For the hard-core model the well-known works of Kahn~\cite{Kah01}, Zhao~\cite{Zha10}, and Galvin--Tetali~\cite{GT04} use the entropy method and the `bipartite swapping trick' to show that over $\Delta$-regular graphs the complete bipartite graph $K_{\Delta,\Delta}$ maximizes the free energy density for the entire parameter range. 
This was strengthened to the level of occupancy fraction in~\cite{DJPR17a}.
A result of Cutler and Radcliffe~\cite{CR14} gives that $K_{\Delta+1}$ is the minimizer of both the free energy and the occupancy fraction (see~\cite{DJPR18b}).

The relevant free energy maximization result for the antiferromagnetic Ising model is~\cite[Corollary 1.15]{SSSZ20}, showing that $K_{\Delta,\Delta}$ is indeed the maximum in the entire parameter range. 
The key to interpreting their result in our setting is the observation that $Z_G(B,\lambda)$ is a weighted sum over graph homomorphisms.
Our Theorem~\ref{thm:d3max} strengthens their result but only in a subset of the parameter range.
We conjecture the extension to the entire range and to all $\Delta$ in Conjecture~\ref{conj:d3}\ref{itm:conjd3max}.
There are several other instances of models where the extremes of the free energy density or occupancy fraction over regular graphs are known or conjectured. See~\cite{CCPT17,Csi17a,Zha17}.

\section{Proof sketch}

The occupancy method starts by showing that the minimum (or maximum) occupancy fraction over a class of graphs is the solution of a constrained optimization problem whose variables are defined by spin systems on graphs. 
Typically, one chooses a spin assignment $\sigma$ distributed according to the model and randomly chooses a small subgraph $F$ of the graph, and then reveals the $\sigma$ restricted to some prescribed subset $S\subset V(F)$. 
We call the outcome of this experiment a \emph{local view}.
The variables of the optimization problem are then the probabilities of the  local views in the distribution corresponding to this experiment. 
If one does this carefully then the occupancy fraction is a linear function of these variables.
In this setup, minimizing the occupancy fraction over a class of graphs corresponds to finding the minimum occupancy fraction over the feasible set of variables which arise when one considers the spin model on graphs in the class. 
We do not know of general methods for solving such optimization problems as the feasible set is potentially challenging to express, but one can often make progress by finding \emph{linear} constraints satisfied on the entire feasible set and solving the linear programming relaxation of the graph-defined optimization problem. 
This expands the feasible set, but if the optimum over the enlarged feasible region can be shown to arise from the spin model on a graph in the desired class, then a solution of the relaxation must be a solution of the original problem.

Simple examples of this proof idea appear in~\cite{CPT17,DdKP21,DJPR17a,DJPR17b,DJPR18a,DKPS20,Zha17}. In the most basic cases the relaxed optimization problems are amenable to elementary analysis and do not require serious computation. 
For more complex examples, e.g.~\cite{Dav18,PP18}, the sheer number of variables becomes unmanageable and one turns to a computer analysis. 
In fact, merely enumerating all the required variables is a nontrivial combinatorial problem in its own right, and as in the case of~\cite{Dav18} we turn to advanced `canonical isomorph-free generation' techniques~\cite{McK98} to attack the problem.

Aside from the practical issues of scaling one's approach to large numbers of variables, there are two main problems that can arise. Firstly, while any given finite linear program can be solved in finite time, we have an infinite family of linear programs parameterized by the relevant variables $\Delta$, $B$ and $\lambda$. 
It is not necessarily straightforward to establish for the entire family that some candidate solution (as a function of $\Delta$, $B$, $\lambda$) is indeed optimal. 
The standard approach is to exploit LP duality and construct a feasible solution for the dual whose objective value is the occupancy fraction of the conjectured graph, and to show that for all parameters in some range that desired solution is dual feasible.
This typically involves solving a linear system defined by the active primal constraints to find dual variables, and proving that a large number of inequalities involving polynomials in the parameters $B$ and $\lambda$ all hold.
A second, more fundamental problem is that there can be a gap between the true graph-defined optimal value and the optimal value of the linear relaxation. 
We confront each of these problems in our study of the antiferromagnetic Ising model and we are only partially successful in solving these issues, which is why our main result applies only to a subset of the desired parameter space.
A novel difficulty in our case is that exact computation requires handling algebraic reals. 
In previous works (e.g.\ \cite{DJPR17a,Dav18,PP18}) one can choose rational values of the parameters and employ exact rational LP solvers to investigate the problem (e.g.\ to find active primal constraints), but we are forced to work with algebraic reals as $\lam_c(\Delta,B)$ is not necessarily rational when $\Delta\ge3$ is an integer and $B\in[0,1]\cap\mathbb{Q}$. 
This necessitates the use of more flexible combinatorial optimization software and costs significantly more computation time.

\section{Proof}

In this section we prove Theorems~\ref{thm:d3lc},~\ref{thm:d3min} and~\ref{thm:d3max} with the assistance of a computer.
The random variable we study is given by the following experiment.
Fix $\Delta$, $B$, $\lambda$ and let $G$ be a $\Delta$-regular graph. 
Select a vertex $\mathbf u$ of $G$ uniformly at random and a spin assignment $\bm\sigma$ from the Ising model $\mu_{G,B,\lambda}$. 
Then let the random variable $\mathbf L$ be the induced subgraph $G[N^2[\mathbf u]]$ together with the spin assignment $\bm \tau = \bm\sigma|_{N^2(\mathbf u)}$ restricted to the vertices at distance 2 from $\mathbf u$. 
It is important to retain the information that $\mathbf u$ in $\mathbf L$ was the vertex sampled from the graph, so we may consider $\mathbf L$ a rooted graph with $\mathbf u$ identified as the root, as well as having a spin assignment for $N^2(u)$. 
We call such a rooted graph with spin assignment a \emph{local view} as it provides a local view of some graph structure of $G$ and the spin assignment $\sigma$. 
Formally, a local view is then a tuple $L=(H,u,\tau)$ where $H$ is a graph containing the specified vertex $u$, and $\tau$ is a spin assignment to a (possibly empty) subset of the vertices of $H$. 

In all graphs $H$ of interest, the specified $u$ will have $\Delta$ neighbors, and each of these neighbors also has $\Delta$ neighbors, but beyond that we must enumerate all possible rooted graphs $(H,u)$ and spin assignments $\tau$ to $N^2(u)$ in $H$.
In fact, for the Ising model on 3-regular graphs there are exactly 23 local views up to isomorphism (where we consider only graph isomorphisms that fix the root). 

Each $\Delta$-regular graph $G$ gives us a probability distribution $x=x(G,B,\lambda)$ over $\mathcal L$ which we represent as a vector in $[0,1]^{\mathcal L}$. 
As in the proof sketch, we establish a set of linear constraints on the feasible vectors $x$. 
Our constraints arise from the fact that in any regular graph a uniform random neighbor $\mathbf v$ of a uniform random vertex $\mathbf u$ is itself distributed uniformly over the vertices.
Thus, we have two ways of computing statistics of a uniform random vertex given the distribution $x$ of the random local view. Our constraints arise from equating the two ways of writing these statistics. 
For convenience, we define for a local view $L=(H,u,\tau)$ the distribution $\mu_L$ on spin assignments of $H$ as the Ising model $\mu_{H,B,\lambda}$ on $H$ conditioned on the spins of $N^2(u)$ agreeing with $\tau$. 
We deliberately suppress the dependence of $\mu_L$ on $B$ and $\lambda$ as they remain fixed throughout. 

One can express the occupancy fraction $\alpha_G(B,\lambda)$ as the probability that a uniform random vertex $\mathbf u$ receives spin $+$ under $\mu_{G,B,\lam}$. 
By the spatial Markov property of the Ising model, this is equal to the probability that in the random local view $\mathbf L=(H,u,\tau)$ distributed according to $x(G,B,\lambda)$, the vertex $u$ gets spin $+$ under the measure $\mu_L$ on spin assignments to $H$.
That is, we define $\alpha_L$ to be the probability that under $\mu_L$, the vertex $u$ in the local view $L$ gets spin $+$ and we have 
\[ \alpha_G(B,\lambda) = \sum_{L\in\mathcal L}\alpha_L x_L,\]
where $x$ is the distribution of the random local view $\mathbf L$ in $G$.
Importantly, this gives us a linear objective function which we wish to minimize $\sum_{L\in\mathcal L}\alpha_L x_L$ over feasible $x$.

For linear constraints on the feasible $x$, note that one can also compute the probability that a uniform random neighbor of $u$ receives spin $+$ under $\mu_L$, which is a different computation over the local views. 
This provides one linear equality constraint on the distribution $x$ that we wish to use, though we do not explicitly add this constraint to our linear program because we turn to richer statistics to obtain more constraints. 
For a local view $L=(H,u,\tau)$ and $j\in \{0,1,\dotsc, \Delta\}$, let $\gamma^u(j)_L$ be the probability under $\mu_L$, exactly $j$ neighbors of $u$ get spin $+$. 
Similarly, let $\gamma^{N(u)}(j)_L$ be the expected number of $+$-neighbors of a uniform random neighbor of $u$ under the measure $\mu_L$ on spin assignments. 
Then for any $x$ which arises from running the local view experiment on a graph and any $0\le j\le \Delta$,
\[ \sum_{L\in\mathcal L}\left(\gamma^u(j)_L -\gamma^{N(u)}(j)_L\right) x_L = 0, \]
giving $\Delta+1$ linear constraints on $x$. 
Constraints of this type first appear in~\cite{DJPR17a}.
Together with the facts that the entries of $x$ are nonnegative and sum to one, we can express a linear programming relaxation of the extremal problem of minimizing the occupancy fraction.

\begin{definition}
    For a given set $\mathcal L$ of local views and set $J$ of constraint indices, the \emph{primal occmin program} is the linear program with variable $x\in\mathbb R^{\mathcal L}$ given by
\begin{align*}
    \min\;&\sum_{L\in\mathcal L}\alpha_L x_L \text{ s.t.} 
    \\&\sum_{L\in\mathcal L} x_L = 1
    \\&\sum_{L\in\mathcal L}\left(\gamma^u(j)_L -\gamma^{N(u)}(j)_L\right) x_L = 0 \qquad \forall j\in J
    \\&x_L \ge 0  \qquad \forall L\in \mathcal L.
\end{align*}
The \emph{primal occmax program} is identical except for the fact that the objective is to be maximized.

Given $\mathcal L$ and $J$ as above, the \emph{dual occmin program} is the linear program with variable $y\in\mathbb R^{\{p\} \cup J}$ given by
\begin{align*}
    \max\;& y_p\text{ s.t.} 
    \\& y_p + \sum_{j\in J}\left(\gamma^u(j)_L -\gamma^{N(u)}(j)_L\right) y_j \le \alpha_L \qquad \forall L\in\mathcal L.
\end{align*}
Note that we use $p$ to index the dual variable corresponding to the primal constraint $\sum_{L\in\mathcal L}x_L=1$ and the set $J$ to index the other dual variables.
We define the corresponding \emph{dual occmax program} analogously.
\end{definition}

\subsection{Minimization: Theorem~\ref{thm:d3lc} and Theorem~\ref{thm:d3min}}
Since $K_4$ has diameter 1, one of our local views is a rooted copy of $K_4$ together with an empty spin assignment (as there are no vertices at distance two from the root), and the vector $x$ corresponding to $K_4$ has a 1 in the entry indexed by $K_4$ and zeros elsewhere. 
Thus, to argue that this $x$ from $K_4$ is optimal in the primal occmin program with local views $\mathcal L$ and $J=\{0,1\}$, it suffices to show that the optimum value achieved in the program with local views $\mathcal L' := \mathcal L\setminus \{K_4\}$ (and the same $J$) is at least the occupancy fraction of $K_4$. 
Given this, in any 3-regular graph without a component isomorphic to $K_4$ we know that the occupancy fraction is at least that of $K_4$. 
If the graph does have components isomorphic to $K_4$ then we can use the fact that for a disjoint union $G=G_1\sqcup G_2$ we have 
\[ |V(G)|\alpha_{G}(B,\lambda) =  |V(G_1)|\alpha_{G_1}(B,\lambda) + |V(G_2)|\alpha_{G_1}(B,\lambda).\]

That is, to prove Theorems~\ref{thm:d3lc} and the statement about $K_4$ being the minimizer in Theorem~\ref{thm:d3min}, it suffices to choose $J=\{0,1\}$ and perform the following steps with the primal and dual occupancy programs with local view set $\mathcal L'$.
\begin{enumerate}
    \item Identify the local views indexing three tight dual constraints, say $\{L_1,L_2,L_3\}\subset \mathcal L'$.
    \item Solve for $y=(y_p,y_0,y_1)$ in the linear system obtained by setting the dual constraints indexed by $(L_1,L_2,L_3)$ to equality.
    \item Check that every dual constraint holds for $y=(\alpha_{K_4}(B,\lambda),y_0,y_1)$.
\end{enumerate}
Note that for this to have any chance of working, we require that the $y_p$ obtained is at least $\alpha_{K_4}(B,\lambda)$. 
This will always be the case when $(B,\lambda)$ are such that indeed $K_4$ minimizes the occupancy fraction over $3$-regular graphs, but even when this does hold our linear program can be a loose relaxation of the graph-theoretic problem and the method fails.

Having done the above steps, we have shown that there is some $y$ feasible in the dual which achieves objective value equal to $\alpha_{K_4}(B,\lambda)$, and hence by LP duality the optimum of the primal is at least $\alpha_{K_4}(B,\lambda)$. 
The first step looks a bit mysterious, but can be done empirically by solving the primal program and looking at the nonnegative values in optimum $x$. 
The reason that we solve for $y=(y_p,y_0,y_1)$ in the three chosen tight dual constraints and then set $y_p$ to $\alpha_{K_4}(B,\lambda)$ is that having taken out the local view $K_4$ the optimum value in the LPs with local view set $\mathcal L'$ could be strictly greater than $\alpha_{K_4}(B,\lambda)$. Thus, setting $y_p$ to the slightly smaller $\alpha_{K_4}(B,\lambda)$ gives us slightly more room to establish the remaining dual constraints. 
While one could imagine other methods for finding dual variables inspired by complementary slackness, this method proved reasonably successful for us.

When following the above steps, the dual constraints (each of which is indexed by some $L\in \mathcal L'$) become inequalities equivalent to 
\[ 0\le \alpha_L - \alpha_{K_4(B,\lambda)} -\sum_{j\in\{0,1\}}\left(\gamma^u(j)_L -\gamma^{N(u)}(j)_L\right) y_j, \]
where each term is a rational function of $B$ and $\lambda$. Thus, provided one can establish that a set of 22 inequalities (19 if you discount the three that we know must hold by construction) of rational functions hold, one can prove Theorem~\ref{thm:d3lc} and the minimziation statement in Theorem~\ref{thm:d3min} for regions of the parameter space.

We use the general-purpose mathematical computing programs Sagemath~\cite{sage} and Wolfram Engine~\cite{wolfram} to assist with the following.
\begin{enumerate}
    \item Generate $\mathcal L'$.
    \item Manually identify a collection $\mathcal T$ of sets of three tight dual constraints.
    \item For each set $T \in \mathcal T$, solve for the relevant dual variables and manually identify a region of the parameter space $(B,\lambda)\in [0,1]^2$ in which we wish to prove symbolically that the dual constraints hold with these variables.
    \item Prove symbolically that in the specified region the dual constraints indeed hold. 
\end{enumerate}
This constitutes a computer-assisted proof of the statement about regions where $K_4$ is the minimizer in Theorem~\ref{thm:d3min}, and Theorem~\ref{thm:d3lc} follows also because the curve given by $\lambda = \lam_c(\Delta,B)$ lies in the union of the regions in which we show dual feasibility for $B\in (0,0.3128)$.

While we could in principle continue for larger $B$, at some point shortly after $B=0.3128$ the optimum value in the LP is smaller than that of $K_4$. We verified this for $B=0.32$ and $\lam_c(3,B)$ over the algebraic reals symbolically with Sagemath.

The statement about $K_4$ minimizing the free energy density follows immediately from the occupancy fraction result and the fact that the union of our regions is downward-closed in the sense that for any $(B,\lambda)$ in the union of the regions, it is true that $(B,\ell)$ is in the union for $0\le \ell\le \lambda$.

To complete the proof of Theorem~\ref{thm:d3min}, we use a computer to search through all $3$-regular graphs on at most 14 vertices and plot the obtained minimziers. 
While 14 vertices is not necessarily representative of the behavior one could observe in general, it is already enough to prove that the minimizer is not always $K_4$, and also that it cannot always $K_4$ or the Petersen graph. 

\subsection{Maximization: Theorem~\ref{thm:d3max}}

\begin{figure}[ht]
    \vspace*{2ex}
    \includegraphics[width=0.2\textwidth]{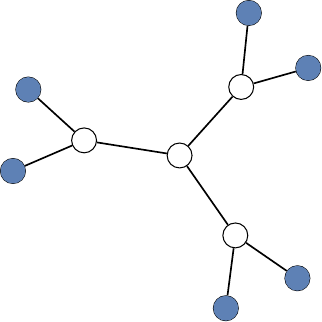}
    \hspace*{0.1\textwidth}
    \includegraphics[width=0.2\textwidth]{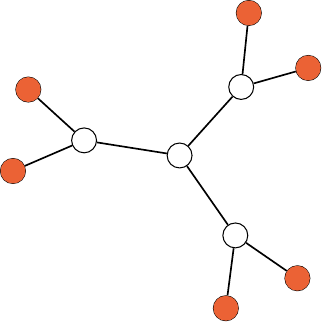}
    \hspace*{0.1\textwidth}
    \includegraphics[width=0.2\textwidth]{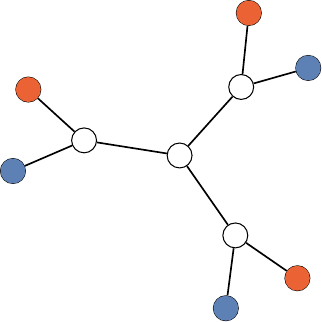}
    \\\vspace*{2ex}
\caption{Local views one can obtain in $K_{3,3}$. The central vertex is the root $u$ of the local view, and vertices assigned $+$ or $-$ are colored red or blue respectively.}
\label{fig:K33lvs}
\end{figure}

The proof of Theorem~\ref{thm:d3max} is analogous. 
This time, it is not the case that a single local view represents $K_{3,3}$, but it is true that disjoint unions of $K_{3,3}$ are the only $3$-regular graphs in which the three neighbors of the root $u$ of the local view always see the same spins on the vertices at distance two from $u$. That is, the only three local views one can obtain in $K_{3,3}$ are the ones shown in Figure~\ref{fig:K33lvs}, and disjoint unions of $K_{3,3}$s are the only graphs with this property. 
One can check that the only feasible solution to the LP supported on these three local views is indeed the exact distribution on local views one gets in $K_{3,3}$. That is, our constraints are enough to capture this key property.

Then we consider only the set of tight constraints corresponding to those three local views, and perform the same kind of dual feasibility analysis over the entire local view set $\mathcal L$. 
The only differences are that the $y_p$ obtained is $\alpha_{K_{3,3}(B,\lambda)}$ and that the dual constraint inequalities are in the opposite direction. That is, we must check inequalities of the form 
\[ 0\ge \alpha_L - \alpha_{K_{3,3}(B,\lambda)} -\sum_{j\in\{0,1\}}\left(\gamma^u(j)_L -\gamma^{N(u)}(j)_L\right) y_j \]
to verify dual feasibility.

\section{Discussion and extensions of the method}

There are several ways one could hope to push our methods further, and we identify a few of the ideas we investigated below. 

In principle, one could consider every subset $T\in{\mathcal L',3}$, solve those dual constraints for equality, and find the maximal region of $[0,1]^2$ in which the corresponding deal feasibility inequalities hold. 
We decided not to perform this rather extensive computation as we could prove that this does not cover the entire parameter space, leaving out for example $B=0.32$ and $\lam=\lam_c(3,0.32)$. 

We solved numerically for the four maximal regions the sets of tight constraints we considered in order to plot Figure~\ref{fig:d3posmin}, but we did not successfully prove symbolically that dual feasibility holds in these exact regions, settling for a slightly smaller region bounded by simpler functions of $B$ and $\lambda$. 
We found that Wolfram Engine could show dual feasibility rather quickly in regions bounded by linear or quadratic functions and so manually fit such regions inside the plotted ones. 
It appears that more specialized tools and significant computing power would be required to prove symbolically that the minimizer is $K_4$ for the entire region in which the linear program is a tight relaxation of the graph-defined optimization problem.

The setup for occupancy fraction optimization that we employ is rather flexible, and one can define `deeper' local views that capture more graph-theoretic structure in the problem. 
We applied the proof technique to local views of depth three in $3$-regular graphs and found (after significant engineering effort and computation time) that there are 2637 local views in this case. 
This is just small enough that one can investigate the corresponding primal occmin program numerically, but the overhead of exact rational or algebraic real computation proved too much to handle.
While we believe that this program could enlarge the known region in $[0,1]^2$ where $K_4$ is optimal, it seems computationally infeasible to extend our dual feasibility analysis to this many local views. 
Not only would there be roughly 100 times more inequalities to check for dual feasibility, each such inequality is now a much larger polynomial with significantly higher degree than in the case of local views of depth 2. 
Initial numerical experiments suggest that the deeper program can show $K_4$ to be the minimizer in a larger region, but that one cannot expect this program to be a tight relaxation in the entire parameter range.
For similar reasons of scale, we did not attempt to tackle $4$-regular graphs at any depth.

It would be interesting to use our methods to show that a graph other than $K_4$ is the minimizer at some point in the parameter space. Unfortunately, it seems that the primal occmin program with local views at depth 2 or at depth 3 does not not yield tight relaxation of the graph-defined optimization problem in e.g.\ the region where we know that the Petersen graph has a smaller occupancy fraction.

We conclude with a question inspired by the fact that in the case of the hard-core model one can study approximate counting with and without a constraint on the number of + spins in the class of triangle-free graphs of bounded degree~\cite{DP23}. 
While the precise minimizing graph is not known, an asymptotically tight bound is known at the critical $\lambda$ value~\cite{DJPR17b}

\begin{question}
    What is the minimum occupancy fraction of the antiferromagnetic Ising model in triangle-free cubic graphs?
\end{question}

\bibliographystyle{habbrv}
\bibliography{bib}

\begin{thebibliography}{10}
\expandafter\ifx\csname url\endcsname\relax
  \def\url#1{\texttt{#1}}\fi
\expandafter\ifx\csname doi\endcsname\relax
  \def\doi#1{\burlalt{\textsc{doi}:\detokenize{#1}}{https://dx.doi.org/#1}}\fi
\expandafter\ifx\csname urlprefix\endcsname\relax\def\urlprefix{\textsc{url:}}\fi
\expandafter\ifx\csname href\endcsname\relax
  \def\href#1#2{#2}\fi
\expandafter\ifx\csname burlalt\endcsname\relax
  \def\burlalt#1#2{\href{#2}{#1}}\fi

\bibitem{AL20a}
V.~L. Alev and L.~C. Lau.
\newblock Improved analysis of higher order random walks and applications.
\newblock In {\em Proceedings of the 52nd {{Annual ACM SIGACT Symposium}} on {{Theory}} of {{Computing}}}, pages 1198--1211, Chicago IL USA, June 2020. ACM.
\newblock \doi{10.1145/3357713.3384317}.

\bibitem{BD97}
R.~Bubley and M.~Dyer.
\newblock Path coupling: {{A}} technique for proving rapid mixing in {{Markov}} chains.
\newblock In {\em Proceedings 38th {{Annual Symposium}} on {{Foundations}} of {{Computer Science}}}, pages 223--231, Miami Beach, FL, USA, 1997. IEEE Comput. Soc.
\newblock \doi{10.1109/SFCS.1997.646111}.

\bibitem{CDKP22}
C.~Carlson, E.~Davies, A.~Kolla, and W.~Perkins.
\newblock Computational thresholds for the fixed-magnetization {{Ising}} model.
\newblock In {\em Proceedings of the 54th {{Annual ACM SIGACT Symposium}} on {{Theory}} of {{Computing}}}, pages 1459--1472, Rome Italy, June 2022. ACM.
\newblock \doi{10.1145/3519935.3520003}.

\bibitem{CCPT17}
E.~Cohen, P.~Csikv{\'a}ri, W.~Perkins, and P.~Tetali.
\newblock The {{Widom}}--{{Rowlinson}} model, the hard-core model and the extremality of the complete graph.
\newblock {\em European Journal of Combinatorics}, 62:70--76, May 2017.
\newblock \doi{10.1016/j.ejc.2016.11.003}.

\bibitem{CPT17}
E.~Cohen, W.~Perkins, and P.~Tetali.
\newblock On the {{Widom}}--{{Rowlinson Occupancy Fraction}} in {{Regular Graphs}}.
\newblock {\em Combinatorics, Probability and Computing}, 26(2):183--194, Mar. 2017.
\newblock \doi{10.1017/S0963548316000249}.

\bibitem{Csi17a}
P.~Csikv{\'a}ri.
\newblock Extremal regular graphs: The case of the infinite regular tree.
\newblock {\em arXiv:1612.01295 [math]}, May 2017, \burlalt{arXiv:1612.01295}{https://arxiv.org/abs/1612.01295}.

\bibitem{CR14}
J.~Cutler and A.~Radcliffe.
\newblock The maximum number of complete subgraphs in a graph with given maximum degree.
\newblock {\em Journal of Combinatorial Theory, Series B}, 104:60--71, Jan. 2014.
\newblock \doi{10.1016/j.jctb.2013.10.003}.

\bibitem{Dav18}
E.~Davies.
\newblock Counting {{Proper Colourings}} in 4-{{Regular Graphs}} via the {{Potts Model}}.
\newblock {\em The Electronic Journal of Combinatorics}, 25(4):P4.7, Oct. 2018.
\newblock \doi{10.37236/7743}.

\bibitem{DdKP21}
E.~Davies, R.~{de Joannis de Verclos}, R.~J. Kang, and F.~Pirot.
\newblock Occupancy fraction, fractional colouring, and triangle fraction.
\newblock {\em Journal of Graph Theory}, 97(4):557--568, 2021.
\newblock \doi{10.1002/jgt.22671}.

\bibitem{DJPR17a}
E.~Davies, M.~Jenssen, W.~Perkins, and B.~Roberts.
\newblock Independent sets, matchings, and occupancy fractions.
\newblock {\em Journal of the London Mathematical Society}, 96(1):47--66, Aug. 2017.
\newblock \doi{10.1112/jlms.12056}.

\bibitem{DJPR17b}
E.~Davies, M.~Jenssen, W.~Perkins, and B.~Roberts.
\newblock On the average size of independent sets in triangle-free graphs.
\newblock {\em Proceedings of the American Mathematical Society}, 146(1):111--124, July 2017.
\newblock \doi{10.1090/proc/13728}.

\bibitem{DJPR18a}
E.~Davies, M.~Jenssen, W.~Perkins, and B.~Roberts.
\newblock Extremes of the internal energy of the {{Potts}} model on cubic graphs.
\newblock {\em Random Structures \& Algorithms}, 53(1):59--75, Aug. 2018.
\newblock \doi{10.1002/rsa.20767}.

\bibitem{DJPR18b}
E.~Davies, M.~Jenssen, W.~Perkins, and B.~Roberts.
\newblock Tight bounds on the coefficients of partition functions via stability.
\newblock {\em Journal of Combinatorial Theory, Series A}, 160:1--30, Nov. 2018.
\newblock \doi{10.1016/j.jcta.2018.06.005}.

\bibitem{DKPS20}
E.~Davies, R.~J. Kang, F.~Pirot, and J.-S. Sereni.
\newblock Graph structure via local occupancy.
\newblock {\em arXiv:2003.14361 [math]}, Mar. 2020, \burlalt{arXiv:2003.14361}{https://arxiv.org/abs/2003.14361}.

\bibitem{DP23}
E.~Davies and W.~Perkins.
\newblock Approximately {{Counting Independent Sets}} of a {{Given Size}} in {{Bounded-Degree Graphs}}.
\newblock {\em SIAM Journal on Computing}, 52(2):618--640, Apr. 2023.
\newblock \doi{10.1137/21M1466220}.

\bibitem{FV17}
S.~Friedli and Y.~Velenik.
\newblock {\em Statistical {{Mechanics}} of {{Lattice Systems}}: {{A Concrete Mathematical Introduction}}}.
\newblock Cambridge University Press, 1 edition, Nov. 2017.
\newblock \doi{10.1017/9781316882603}.

\bibitem{GT04}
D.~Galvin and P.~Tetali.
\newblock On weighted graph homomorphisms.
\newblock In {\em Graphs, Morphisms and Statistical Physics}, volume~63 of {\em {{DIMACS Ser}}. {{Discrete Math}}. {{Theoret}}. {{Comput}}. {{Sci}}.}, pages 97--104. Amer. Math. Soc., Providence, RI, 2004.
\newblock \doi{10.1090/dimacs/063/07}.

\bibitem{HK24}
D.~G. Harris and V.~Kolmogorov.
\newblock Parameter estimation for {{Gibbs}} distributions.
\newblock {\em ACM Trans. Algorithms}, July 2024.
\newblock \doi{10.1145/3685676}.

\bibitem{JMPV23}
V.~Jain, M.~Michelen, H.~T. Pham, and T.-D. Vuong.
\newblock Optimal mixing of the down-up walk on independent sets of a given size.
\newblock In {\em 2023 {{IEEE}} 64th {{Annual Symposium}} on {{Foundations}} of {{Computer Science}} ({{FOCS}})}, pages 1665--1681, Nov. 2023.
\newblock \doi{10.1109/FOCS57990.2023.00101}.

\bibitem{JPSS22}
V.~Jain, W.~Perkins, A.~Sah, and M.~Sawhney.
\newblock Approximate counting and sampling via local central limit theorems.
\newblock In {\em Proceedings of the 54th {{Annual ACM SIGACT Symposium}} on {{Theory}} of {{Computing}}}, {{STOC}} 2022, pages 1473--1486, New York, NY, USA, June 2022. Association for Computing Machinery.
\newblock \doi{10.1145/3519935.3519957}.

\bibitem{Kah01}
J.~Kahn.
\newblock An {{Entropy Approach}} to the {{Hard-Core Model}} on {{Bipartite Graphs}}.
\newblock {\em Combinatorics, Probability and Computing}, 10(3):219--237, May 2001.
\newblock \doi{10.1017/S0963548301004631}.

\bibitem{KPPY24}
A.~Kuchukova, M.~Pappik, W.~Perkins, and C.~Yap.
\newblock Fast and {{Slow Mixing}} of the {{Kawasaki Dynamics}} on {{Bounded-Degree Graphs}}.
\newblock In {\em Approximation, {{Randomization}}, and {{Combinatorial Optimization}}. {{Algorithms}} and {{Techniques}} ({{APPROX}}/{{RANDOM}} 2024)}, pages 56:1--56:24. Schloss Dagstuhl -- Leibniz-Zentrum f{\"u}r Informatik, 2024.
\newblock \doi{10.4230/LIPIcs.APPROX/RANDOM.2024.56}.

\bibitem{McK98}
B.~D. McKay.
\newblock Isomorph-{{Free Exhaustive Generation}}.
\newblock {\em Journal of Algorithms}, 26(2):306--324, Feb. 1998.
\newblock \doi{10.1006/jagm.1997.0898}.

\bibitem{PP18}
G.~Perarnau and W.~Perkins.
\newblock Counting independent sets in cubic graphs of given girth.
\newblock {\em Journal of Combinatorial Theory, Series B}, 133:211--242, Nov. 2018.
\newblock \doi{10.1016/j.jctb.2018.04.009}.

\bibitem{SSSZ20}
A.~Sah, M.~Sawhney, D.~Stoner, and Y.~Zhao.
\newblock A reverse {{Sidorenko}} inequality.
\newblock {\em Inventiones mathematicae}, Mar. 2020.
\newblock \doi{10.1007/s00222-020-00956-9}.

\bibitem{sage}
W.~A. Stein~et al.
\newblock {S}age {M}athematics {S}oftware ({V}ersion 10.4), 2024.
\newblock \urlprefix\url{www.sagemath.org}.

\bibitem{wolfram}
{Wolfram Research, Inc.}
\newblock Wolfram {E}ngine ({V}ersion 14.1), 2024.
\newblock \urlprefix\url{www.wolfram.com/engine}.

\bibitem{Zha10}
Y.~Zhao.
\newblock The {{Number}} of {{Independent Sets}} in a {{Regular Graph}}.
\newblock {\em Combinatorics, Probability and Computing}, 19(2):315--320, Mar. 2010.
\newblock \doi{10.1017/S0963548309990538}.

\bibitem{Zha17}
Y.~Zhao.
\newblock Extremal {{Regular Graphs}}: {{Independent Sets}} and {{Graph Homomorphisms}}.
\newblock {\em The American Mathematical Monthly}, 124(9):827, 2017.
\newblock \doi{10.4169/amer.math.monthly.124.9.827}.

\end{thebibliography}

\appendix
\section{Running our code}

To verify our claims with a computer you will need Sagemath, see \url{sagemath.org}. We tested versions 10.2, 10.3, and 10.4. 
You will also need Wolfram Engine, see \url{www.wolfram.com/engine}. We tested versions 14.0 and 14.1. 
Our code is available at \url{github.com/ed359/IsingOccupancy/tree/d3paper} and contains relevant instructions.

\end{document}